\documentclass{article}

\usepackage[utf8]{inputenc}
\usepackage{algorithm}
\usepackage{algorithmic}
\usepackage{floatflt}
\usepackage{graphicx}
\usepackage{subfigure} 
\usepackage{indentfirst} 

\usepackage{url}
\usepackage{fancyhdr}
\usepackage{amsmath}
\usepackage{amsfonts}
\usepackage{amssymb}		
\usepackage{amsthm}
\usepackage{float}
\usepackage{latexsym}
\usepackage{graphics}
\usepackage{verbatim}
\usepackage[T1]{fontenc}
\usepackage{lmodern}
\usepackage{enumerate}
\usepackage{setspace}
\usepackage{tikz-cd}
\usepackage{multirow}
\usepackage{geometry}

\usepackage[width=.8\textwidth, format=hang]{caption} 

\usepackage{color}

\newcommand{\blue}[1]{{\color{blue} #1}}

\newcommand{\FF}{{\mathbb{F}}}

\newcommand{\ZZ}{{\mathbb{Z}}}

\newcommand{\ord}{\text{ord}}

\theoremstyle{definition}
\newtheorem{theorem}{\bf \blue{Theorem}}
\newtheorem{lemma}[theorem]{\bf \blue{Lemma}}

\newtheorem{remark}[theorem]{\bf \blue{Remark}}

\newtheorem{problem}[theorem]{\bf \blue{Problem}}

\providecommand{\keywords}[1]
{
  \textbf{\textit{Keywords---}} #1\\
}

\providecommand{\msc}[1]
{
  \textbf{\textit{MSC2020---}} #1
}

\title{Root Extraction in Certain Finite Abelian $p$-Groups\footnote{This work appeared as a chapter in the author's Ph.D. thesis, titled \emph{Isogeny-based Quantum Resistant Undeniable Blind Signature and Authenticated Encryption Schemes}, 2018.}}

\author{M. S. Srinath\\srinathms@sssihl.edu.in\\Department of Mathematics and Computer Science\\Sri Sathya Sai Institute of Higher Learning, Prasanthi Nilayam\\Puttaparthi, 515134, Andhra Pradesh, India}
\date{}
\begin{document}

\maketitle


\begin{abstract}
  We formulate a problem called \emph{Generalized Root Extraction} in finite Abelian groups that have more than one generator. We then study this problem for the specific case of the torsion subgroups of elliptic curves. We give a necessary and sufficient condition for the existence of a solution. We also present an algorithm to find a solution. Our algorithm easily generalizes to Abelian groups of prime power order that have a specific structure. We then discuss a variant of this problem called \emph{Simultaneous Root Extraction} and present an algorithm for solving it.
\end{abstract}

\keywords{Finite Abelian $p$-Groups, Root Extraction, Simultaneous Root Extraction.}

\msc{20-08, 20K01, 14H52.}


\section{Introduction}

Let $G$ be a finite Abelian group. Given $h \in G$ and $k$ a positive integer, the problem of \emph{root extraction} in multiplicative notation is to find $g$ such that 
\[ h = g^k. \]
In this case, we say that $g$ is $k$-th root of $h$. Maurer and Wolf~\cite[Theorem~11]{maurer} have given an algorithm for extracting roots in cyclic groups. Damg{\aa}rd and Koprowski~\cite{Damgard2002} have given generic lower bounds for this problem in generic groups.

Root extraction in finite fields has been studied by Cao et al.~\cite{Cao2012} and Koo et al.~\cite{Koo2013} as an improvement over Adleman-Manders-Miller algorithm~\cite{AMM}. Root extraction in matrix rings over fields has been studied by Otero~\cite{Otero1990}. In this work, we study root extraction in Abelian groups of prime power order that are of the form $G \approx \prod_{i=1}^{N}{\frac{\ZZ}{p^{e}\ZZ}}$.

The structure theorem for finite Abelian groups \cite[Theorem 11.1]{gallian} states that every finite Abelian group is a direct product of cyclic groups of prime-power order. More precisely, let $G$ be a finite Abelian group. Then, there exists positive integers $e_i$ and primes $p_i$ (not necessarily distinct) for $i=1, \ldots, N$ such that
\begin{equation}
G \approx \prod_{i=1}^{N}{\frac{\ZZ}{p_i^{e_i}\ZZ}}
\end{equation}
and that this decomposition is unique (up to permutations of $p_i$).

Hence, it seems natural to ask the question, given $h \in G$ and positive integers $k_1, k_2, \ldots, k_N$, find generators $g_1, g_2, \ldots, g_N$ such that 
\begin{equation}
h = g_1^{k_1}g_2^{k_2} \ldots g_N^{k_N}.
\end{equation}
We will call this the \emph{generalized root extraction problem}.

\section{Formulation of the Problem}

We will give a formulation for the specific case of torsion subgroups of elliptic curves. We consider this specific case as they are of interest in isogeny-based cryptography~\cite{feo}. We believe that our formulation and solution throw more light into the understanding of the emerging area of isogeny-based cryptography.

Let $p$ be a prime of the form $\ell^e\cdot f\pm 1$ where $\ell$ is a small prime, $e$ is a positive integer and $f$ is the co-factor. Let $E$ be a supersingular elliptic curve defined over $\FF_{p^2}$ and having cardinality $(\ell^e\cdot f)^2$. The torsion subgroup $E[\ell^e]$ has the structure
\begin{equation}
E[\ell^e] \approx \frac{\ZZ}{\ell^e\ZZ} \times \frac{\ZZ}{\ell^e\ZZ}.
\end{equation}
Let $P,Q \in E[\ell^e]$ be two points of order $\ell^e$ that generate $E[\ell^e]$. Let $K \in E[\ell^e]$ be another point with $K = mP + nQ$ for $m,n \in \ZZ/\ell^e\ZZ$.

\begin{problem}[Generalized Root Extraction Problem, GREP]
    \label{def:grep}
	Given $m, n \in \ZZ/\ell^e\ZZ$ and $K \in E[\ell^e]$; determine $P$ and $Q$ such that $K = mP + nQ$ and $\left<P,Q\right> = E[\ell^e]$.
\end{problem}

For the rest of this work (except Theorem~\ref{thm:root}), we restrict our discussion to the above case, i.e., the torsion subgroups of supersingular curves that are of interest in isogeny-based cryptography. Additionally, we seek solutions that generate the torsion subgroup, keeping in mind the requirements of the possible cryptographic applications. In section~\ref{sec:repalg}, we present an algorithm for solving the generalized root extraction problem. At the end of the section, we also remark on how to generalize this for groups of the form $G \approx \prod_{i=1}^{N}{\frac{\ZZ}{p^{e}\ZZ}}$. In section~\ref{sec:repdiscussion}, we consider simultaneous equations involving roots and then give a solution to this system.

Before discussing the algorithm, we prove a lemma and a theorem. The lemma below is an easy observation while the theorem is a generalization of~\cite[Theorem~11]{maurer}. Both of these results will be used in the algorithm.

\begin{lemma} \label{lem:ord}
	Let $P$ and $Q$ be generators of $E[\ell^e]$. Suppose $K = mP + nQ$ where $\ell \nmid \gcd(m,n)$. Then, $\ord(K) = \ell^e$.
\end{lemma}

\begin{proof}
	Suppose $\ord(K) = \ell^{e_1}$, for some $e_1<e$. Then,
	\begin{align*}
	O = \ell^{e_1}K &= \ell^{e_1}mP + \ell^{e_1}nQ
	\end{align*}
	where $O$ is the identity of the elliptic curve group $E$. Since, $\ord(P) = \ord(Q) = \ell^e$ and $\ell \nmid \gcd(m,n)$, at least one term on the right-hand side of the above equation is non-zero. This implies that $P$ and $Q$ are linearly dependent which is a contradiction as they are generators of $E[\ell^e]$. Therefore, $\ord(K) = \ell^e$.
\end{proof}

We now prove a theorem that generalizes a result due to Maurer and Wolf~\cite[Theorem~11]{maurer}. We then apply this result specifically for the torsion group $E[\ell^e]$.

\begin{theorem}
	\label{thm:root}
	Let $G=\langle g_1,g_2, \dots g_N\rangle$ be an Abelian group such that $|G|=\ell^t\cdot s$ with $\gcd(\ell,s)=1$. For $1\leq r < t$, $\ell^r$-th roots can be computed in $G$ in time $O(L^{N/2}\log |G|)$ where $L$ is the largest prime divisor of $|G|$.
\end{theorem}

\begin{proof}
	Let $h$ be a $\ell^r$-th power in $G$. Using the generalized Pohlig-Hellman algorithm~\cite{teske-ph}, the extended discrete logarithm of $h=\prod_{i=1}^N{g_i^{k_i}}$ can be computed in time $O(L^{N/2}\log |G|)$ where $L$ is the largest prime divisor of $|G|$. Since $h$ is a $\ell^r$-th power in $G$, each $k_i$ $(1 \leq i \leq N)$ is an integer multiple of $\ell^r$. Let $d \equiv -s^{-1} \mod \ell^r$; i.e., $sd + 1 = c\ell^r$ for some integer $c$. Then $x$ is an $\ell^r$-th root of $h$, where
	\[ x = h^{(sd+1)/\ell^r}\cdot\prod_{i=1}^N(g_i^{s\cdot (k_i/\ell^r)\cdot d})^{-1}. \]
	Indeed,
	\begin{align*}
	x^{\ell^r} = {\left[h^{(sd+1)/\ell^r}\cdot\prod_{i=1}^N(g_i^{s\cdot (k_i/\ell^r)\cdot d})^{-1}\right]}^{\ell^r}
	&= h^{(sd+1)}\cdot\prod_{i=1}^N(g_i^{s\cdot (k_i/\ell^r)\cdot d})^{-\ell^r}\\
	&= h^{(sd+1)}\cdot\prod_{i=1}^Ng_i^{-s\cdot k_i\cdot d}\\
	&= h^{c\ell^r}\cdot\prod_{i=1}^Ng_i^{(1-c\ell^r)k_i}\\
	&= h^{c\ell^r}\cdot\prod_{i=1}^Ng_i^{k_i}\cdot\left[\prod_{i=1}^ng_i^{k_i}\right]^{-c\ell^r}\\
	&= h^{c\ell^r}\cdot h\cdot h^{-c\ell^r}\\
	&= h.
	\end{align*}
\end{proof}


\section{The Algorithm}
\label{sec:repalg}

Given a point $K \in E[\ell^e]$ and two integers $m,n \in \frac{\ZZ}{\ell^e\ZZ}$, the aim of the algorithm is to find two points $\{P,Q\}$ such that, $K=mP+nQ$ and $\left<P,Q\right>= E[\ell^e]$.

The first step is to find the order of $K$. Since $E[\ell^e] \approx \ZZ/\ell^e\ZZ \times \ZZ/\ell^e\ZZ$ and $\ell$ is prime, $ \ord(K) \mid \ell^e$. More specifically, $\ord(K)$ is a power of $\ell$. Hence, $\ord(K) = \ell^u$ for some $u \leq e$. Now consider the ordered set, $\{\ell^jK\}$, $0 \leq j \leq e$. The first element in the list will be $K$ and the last element will be $O$, the identity of the elliptic curve group. Further, $\ell^jK \neq O$ for $0 \leq j < u$ and $\ell^jK = O$ for $u \leq j \leq e$. Thus, $\ord(K)$ can be found in at most $e$ steps and hence $O(e)$ steps.

Let $r$ be the highest power of $\ell$ that divides both $m$ and $n$. i.e., $0 \leq r < e$ is an integer such that $\ell^r | \gcd(m,n)$ and $\ell^{r+1} \nmid \gcd(m,n)$. By Lemma~\ref{lem:ord}, GREP is solvable if and only if $u+r=e$.

\begin{theorem}[Existence Theorem]
    \label{thm:existence}
    The solution to GREP described in Problem \ref{def:grep} exists if and only if $u+r=e$.
\end{theorem}

\begin{proof}
    The proof follows immediately from Lemma \ref{lem:ord}.
\end{proof}

We consider two different cases for the algorithm. We present pseudo-codes for both of these cases at the end of the section.

\begin{description}
	\item[Case 1:] $\ell$ does not simultaneously divide both $m$ and $n$ (so $r=0$ and $\ord(K) = \ell^e$).\\ Without loss of generality, assume $\ell \nmid n$.
	\begin{enumerate}
		\item Find another point $K'$ of order $\ell^e$ such that $K$ and $K'$ generate the $E[\ell^e]$ torsion subgroup of $E$.
		An algorithm for doing this is described by Azarderakhsh et al.~\cite[Section 3.2]{keycomp}.
		\item Now assign $P := K'$ and $Q := n^{-1}(K-mK')$.\\ Then $K = mP + nQ$.

As $\ell \nmid n$, we can compute $n^{-1}\mod\ell^e$ efficiently using extended Euclid's algorithm.
	\end{enumerate}
{\bf Claim:} $\{P,Q\}$ generates $E[\ell^e]$.
\begin{proof}
It is easy to check that $K$ and $K'$ are linearly independent, both of order $\ell^e$. By Lemma~\ref{lem:ord}, $\ord(K-mK')=\ell^e$. Therefore, $\ord(Q) =\ell^e $.
Thus, it suffices to show that $K'$ and $K-mK'$ are linearly independent. Consider,
\begin{align*}
aK' + b(K-mK') &= O \\
(a-bm)K' + bK &= O.
\end{align*}
Since, $K, K'$ generate $E[\ell^e]$, $b=0$ and $a - bm = 0$. Hence, $a=0$.\\
Therefore, $\{P, Q\}$ indeed generates $E[\ell^e]$.
\end{proof}
We summarize the algorithm for solving this case in Algorithm~\ref{alg:case1}.

\item[Case 2:] $\ell$ divides both $m$ and $n$.\\ We solve this case by reducing the problem to the previous case. Let $r$ be the highest power of $\ell$ that divides both $m$ and $n$. Now $K$ can be written as\\
$K = \ell^r(m_1P + n_1Q)$ where $\ell^rm_1 = m$, $\ell^rn_1 = n$ and $\ell \nmid \gcd(m_1,n_1)$.\\
Let $R = m_1P + n_1Q$. Therefore, $K=\ell^rR$ and by Lemma~\ref{lem:ord}, $\ord(R) = \ell^e$.

We use Theorem~\ref{thm:root} to find $R$ having known $K$. (Once $R$ is known, we can use Case~1 for finding $P$ and $Q$.)
Since $K=\ell^rR$, $R$ is the $\ell^r$-th root of $K$ in $E[\ell^e]$. For the group $E[\ell^e]$, we let $N=2$, $s=1$, $c=1$ and $d=\ell^r-1$ as per the notation of Theorem~\ref{thm:root}. Choose points $P'$ and $Q'$ that generate $E[\ell^e]$. By substituting the values in Theorem~\ref{thm:root}, we get
\[ R = \left(\frac{k_1}{\ell^r}\right)P' + \left(\frac{k_2}{\ell^r}\right)Q' \]
where $k_1, k_2$ is the solution of the equation $K = k_1P' + k_2Q'$.

Now, solving for the equation $R = m_1P + n_1Q$ in the torsion subgroup $E[\ell^e]$ will yield a solution to the original equation $K = mP + nQ$.

Note that $\ell \nmid \gcd(m_1,n_1)$. Therefore, $\{P,Q\}$ can be calculated by following the steps in Case~1. Thus, we solve Case~2 by using Theorem~\ref{thm:root} and reducing it to Case~1. We summarize the algorithm for solving this case in Algorithm~\ref{alg:case2}.
\end{description}

\begin{algorithm}[!ht]
	\caption{Root Extraction (Case~1)}
	\label{alg:case1}
	Input: $E[\ell^e], K, m, n$.\\
	Precondition: $m,n \in \frac{\ZZ}{\ell^e\ZZ}$, $\ord(K) = \ell^e$ and $\ell \nmid n$.\\
	Output: $P,Q$.\\
	Postcondition: $K=mP+nQ$, $\left<P,Q\right>= E[\ell^e]$.\\
	\textbf{begin}
	\begin{algorithmic}[1]
		\STATE Find $K'$ such that $\left<K,K'\right>= E[\ell^e]$. \hfill $O(2e\log\ell)$
		\STATE Assign $P:=K'$ and $Q := n^{-1}(K-mK')$. \hfill $O(2e\log\ell)$
		\RETURN $P,Q$.
	\end{algorithmic}
\end{algorithm}

\begin{algorithm}[!ht]
	\caption{Root Extraction (Case~2)}
	\label{alg:case2}
	Input: $E[\ell^e], K, m, n$.\\
	Precondition: $m,n \in \frac{\ZZ}{\ell^e\ZZ}$, $K \in E[\ell^e]$, $\ell \mid m$ and $\ell \mid n$.\\
	Output: $P,Q$.\\
	Postcondition: $K=mP+nQ$, $\left<P,Q\right>= E[\ell^e]$.\\
	\textbf{begin}
	\begin{algorithmic}[1]
		\STATE Find $u$ such that $\ord(K) = \ell^u$. \hfill $O(e)$
		\STATE Find $r$ such that $\ell^r | \gcd(m,n)$ and $\ell^{r+1} \nmid \gcd(m,n)$. \hfill $O(e)$
		\IF{$u + r \neq e$} 
			\STATE \textbf{raise exception} \texttt{``ERROR: Solution does not exists.''}
		\ENDIF
		\STATE Assign $m_1 := \frac{m}{\ell^r}$ and $n_1 := \frac{n}{\ell^r}$. 
		\STATE Find $P'$ and $Q'$ such that $\left<P',Q'\right>= E[\ell^e]$. \hfill $O(2e\log\ell)$
		\STATE Find $k_1$ and $k_2$ such that $K = k_1P' + k_2Q'$. \hfill $O(2e\ell\log\ell)$
		\STATE Assign $ R := \left(\frac{k_1}{\ell^r}\right)P' + \left(\frac{k_2}{\ell^r}\right)Q'$. 
		\STATE Run Algorithm~\ref{alg:case1} with the input $E[\ell^e], R, m_1, n_1$. \hfill $O(4e\log\ell)$
		\RETURN $P,Q$.
	\end{algorithmic}
\end{algorithm} 

\begin{remark}[Generalization]
    The algorithm presented above easily generalizes to finite Abelian $p$-groups that are of the form $G \approx \prod_{i=1}^{N}{\frac{\ZZ}{p^{e}\ZZ}}$. Step 7 and Step 8 of Algorithm~\ref{alg:case2} are to be replaced with algorithms by Sutherland~\cite[\S4]{sutherland2011} and Teske~\cite[\S4]{teske-ph} respectively.
\end{remark}


\section{Complexity Analysis}

In the analysis that follows, we count only the number of group operations to obtain the complexity bound.


The analysis of Algorithm~\ref{alg:case1} is fairly simple. In the first step, we find a point $K'$ that is linearly independent to $K$ and has order $\ell^e$. A method for finding such a point is described by Azarderakhsh et al.~\cite[Section 3.2]{keycomp}. The most time-consuming step in their algorithm is computing the Weil pairing of the points in $E[\ell^e]$ which roughly takes $2e\log\ell$ steps and asymptotically has the complexity $O(e\log\ell)$. The second step of our algorithm involves finding $n^{-1}\mod\ell^e$ which can be computed using extended Euclid's algorithm and has complexity $O(e\log\ell)$. Therefore, the overall complexity of Algorithm~\ref{alg:case1} is $O(e\log\ell)$.

In the first step of Algorithm~\ref{alg:case2}, we compute the order of $K$. This can be done in $u$ steps where $\ord(K) = \ell^u$. Since $u$ is bounded above by $e$, the complexity of the first step is $O(e)$.
The second step of the algorithm takes $r$ steps. Since $r\leq e$, the complexity is $O(e)$.
The seventh step involves finding generators of $E[\ell^e]$. A method for finding generators is described by Azarderakhsh et al.~\cite[Section 3.2]{keycomp}. The most time-consuming step in their algorithm is computing the Weil pairing of the points in $E[\ell^e]$ which has the complexity $O(e\log\ell)$. The eighth step involves solving extended discrete logarithm in $E[\ell^e]$. The generalized Pohlig--Hellman algorithm by Teske~\cite{teske-ph} solves it in $O(e\ell\log\ell)$. The tenth and final step involves running Algorithm~\ref{alg:case1} which, by the previous discussion, takes $O(e\log\ell)$ time.

The following theorem summarizes the above discussion on the complexity of our algorithm for solving the GREP in $E[\ell^e]$.

\begin{theorem}
	Let $E/\FF_{p^2}$ be a supersingular elliptic curve with $\#E(\FF_{p^2}) = (\ell^e\cdot f)^2$ where $p$ is a prime of the form $ p = \ell^e\cdot f \pm 1$. Let $K \in E[\ell^e]$ be a point such that $K = mP + nQ$ for some $m,n \in \frac{\ZZ}{\ell^e\ZZ}$ where $\left<P,Q\right> = E[\ell^e]$. Given $m$, $n$ and $K$, the points $P$ and $Q$ can be found in time $O(e\ell\log\ell)$. 
\end{theorem}

\begin{proof}
	Since Algorithm~\ref{alg:case2} calls Algorithm~\ref{alg:case1} as a submodule, the overall complexity of GREP is simply the complexity of Algorithm~\ref{alg:case2}. The complexity of Algorithm~\ref{alg:case2} is asymptotically $O(e) + O(e\log\ell) + O(e\ell\log\ell)$. Since $O(e\ell\log\ell)$ is the dominating term, we are done.
\end{proof}

\begin{remark}[Non-uniqueness of the solution]
    As can be seen from Step\:1 of Algorithm~\ref{alg:case1} and Step\:7 of Algorithm~\ref{alg:case2}, the solution to GREP is not unique.
\end{remark}


\section{Simultaneous Root Extraction}
\label{sec:repdiscussion}
 The non-uniqueness of the solution is due to the dimensionality of the problem. More precisely, we have one equation and two unknowns. In this section, we consider a variant of GREP; namely, simultaneous root extraction of two points. Suppose we need to find the solution to the following system with two equations in two unknowns $P$ and $Q$,
\begin{align}
	K_1 &= m_1P + n_1Q \label{eqn:K1}\\
	K_2 &= m_2P + n_2Q.\label{eqn:K2}
\end{align}

We assume that $m_1n_2 - m_2n_1 \neq 0$ i.e., one equation is not a scalar multiple of the other. By using Theorem~\ref{thm:root}, if necessary rearranging and renaming, we may further assume without loss of generality that $\ell \nmid n_2$. 

We adopt the matrix notation as in~\cite[p.\,194]{siltate}, and write the  equations~\eqref{eqn:K1} and \eqref{eqn:K2} as
\[ (K_1, K_2) = (P, Q) \left( \begin{matrix}
	m_1 & m_2\\
	n_1 & n_2
\end{matrix} \right). \]
\begin{tabular}{r c l}
	The solution $\{P,Q\}$ is unique & $\iff$ & The matrix, $M = \left(\begin{smallmatrix} 
	m_1 & m_2\\
	n_1 & n_2
	\end{smallmatrix} \right)$
	is invertible. \\
	& $\iff$ & $\det M$ is invertible in $\ZZ/\ell^e\ZZ$. \\
	& $\iff$ & $m_1n_2 - m_2n_1$ is a unit in $\ZZ/\ell^e\ZZ$. \\
	& $\iff$ & $\ell \nmid m_1n_2 - m_2n_1$.
\end{tabular}

If $\ell \nmid m_1n_2 - m_2n_1$, then the solution to the simultaneous equations \eqref{eqn:K1} and \eqref{eqn:K2} is unique and is given by
\begin{equation}
(P,Q) = (K_1,K_2)M^{-1}.
\end{equation}

Suppose $\ell \mid m_1n_2 - m_2n_1$, then the solution is not unique. Multiplying the equation~\eqref{eqn:K1} by $n_2$ and the equation~\eqref{eqn:K2} by $n_1$ and subtracting one from the other, we get
\begin{equation} \label{eqn:KK}
	n_2K_1 - n_1K_2 = (m_1n_2 - m_2n_1)P.
\end{equation}

Since $\ell \mid m_1n_2 - m_2n_1$, let $r$ be the highest power of $\ell$ that divides $m_1n_2 - m_2n_1$. Therefore $m_1n_2 - m_2n_1 = \ell^rs$ for some $s \in \ZZ/\ell^e\ZZ$ and $\gcd(\ell,s) = 1$.

Let $H = s^{-1}(n_2K_1 - n_1K_2)$. Now, from equation~\eqref{eqn:KK}, $P$ is an $\ell^r$-th root of $H$ which can be computed in the same manner as in Algorithm~\ref{alg:case2}. Substituting for $P$ in equation~\eqref{eqn:K2}, $Q$ can be efficiently computed as
\[Q = n_2^{-1}(K_2 - m_2P). \]
The points $P$ and $Q$, thus computed, is a solution to the simultaneous equations~\eqref{eqn:K1} and \eqref{eqn:K2} whenever $\ell \mid m_1n_2 - m_2n_1$. We note that the inverses of $s$ and $n_2$ are to be computed modulo $\ell^e$. The algorithm for simultaneous root extraction is summarized in Algorithm~\ref{alg:sre}.

\begin{algorithm}[!ht]
	\caption{Simultaneous Root Extraction}
	\label{alg:sre}
	Input: $E[\ell^e], K_1, m_1, n_1, K_2, m_2, n_2$.\\
	Precondition: $m_1,n_1, m_2,n_2 \in \frac{\ZZ}{\ell^e\ZZ}$; $K_1,K_2 \in E[\ell^e]$; $m_1n_2 - m_2n_1 \neq 0$ and $\ell \nmid n_2$.\\
	Output: $P,Q \in E[\ell^e]$.\\
	Postcondition: $K_1=m_1P+n_1Q$; $K_2=m_2P+n_2Q$ and $\left<P,Q\right>= E[\ell^e]$.\\
	\textbf{begin}
	\begin{algorithmic}[1]
		\STATE Assign $M = \left(\begin{smallmatrix} 
	m_1 & m_2\\
	n_1 & n_2
	\end{smallmatrix} \right)$.
            \IF{$\ell \nmid m_1n_2 - m_2n_1$} 
			\STATE $(P,Q) = (K_1,K_2)M^{-1}$.
            \ELSE
                \STATE Find $s$ such that $m_1n_2 - m_2n_1 = \ell^rs$ and $\gcd(\ell,s) = 1$.
                \STATE Assign $H :=s^{-1}(n_2K_1 - n_1K_2)$.
                \STATE Find $P'$ and $Q'$ such that $\left<P',Q'\right>= E[\ell^e]$. \hfill $O(2e\log\ell)$
		      \STATE Find $k_1$ and $k_2$ such that $H = k_1P' + k_2Q'$. \hfill $O(2e\ell\log\ell)$
		      \STATE Assign $ P := \left(\frac{k_1}{\ell^r}\right)P' + \left(\frac{k_2}{\ell^r}\right)Q'$.
                \STATE Assign $Q := n_2^{-1}(K_2 - m_2P)$.
		\ENDIF
		\RETURN $P,Q$.
	\end{algorithmic}
\end{algorithm}


\section{Acknowledgements}

This work is dedicated to Bhagawan Sri Sathya Sai Baba, the Founder Chancellor of Sri Sathya Sai Institute of Higher Learning, Puttaparthi. The author thanks (i)~Dr.~Uday Kiran, Sri Sathya Sai Institute of Higher Learning, Puttaparthi, (ii)~Prof.~Vijay Patankar, FLAME University, Pune, and (iii)~Dr.~Andrew Sutherland, Massachusetts Institute of Technology (MIT), Greater Boston for their feedback on the earlier drafts of this article.

\bibliographystyle{acm}
\bibliography{References}

\end{document}